\documentclass[11pt]{article}

\usepackage[T2A]{fontenc}
\usepackage[cp1251]{inputenc}
\usepackage[tbtags]{amsmath}
\usepackage{amsfonts,amsthm,amssymb,amstext,amsbsy,
latexsym,mathrsfs,mparhack,marginnote}

\usepackage{mathrsfs}

\usepackage{graphicx}

\usepackage{accents}

\usepackage[dvipsnames]{xcolor} 

\usepackage[bbgreekl]{mathbbol}
\usepackage{bold-extra}

\DeclareSymbolFontAlphabet{\mathbbl}{bbold}

\usepackage{hyperref}
\hypersetup{citecolor=magenta}

\theoremstyle{plain}

\theoremstyle{definition}

\newtheorem*{qroof}{Proof}
\newtheorem*{ackn}{Acknowledgement}
\newtheorem*{Qroblem}  {Problem}

\usepackage{array,longtable}
\newcolumntype{C}{>{$}c<{$}}  
\newcolumntype{L}{>{$}l<{$}}  
\newcolumntype{R}{>{$}r<{$}}  
\numberwithin{theorem}{subsection}

\input{84a.sty}

\begin{document}

\title
{On Petr Novikov's problem of ordered systems  
of uniform sets%
\thanks
{The research was carried out within the state 
assignment of Ministry of Science and Higher 
Education of the Russian Federation for IITP RAS.}
}

\author 
{
Vladimir~Kanovei\thanks{IITP RAS,
  Moscow, Russia, \ 
  {\tt kanovei@iitp.ru} --- contact author. 
}  
\and
Vassily~Lyubetsky\thanks{IITP RAS,
  Moscow, Russia, \ {\tt lyubetsk@iitp.ru} 
}
}

\date 
{\today}

\maketitle

\markboth{В.\,Г.~Кановей и В.\,А.~Любецкий}
{Petr Novikov's problem of ordered systems  
of uniform sets}

\begin{abstract}
We prove that every ordinal $\alpha<\omega_2$ 
is the order type of a certain system of uniform 
Borel sets in the sense of a well-ordering 
relation defined by Petr Novikov. 
This result gives a positive answer to a problem
posed by Nicolas Luzin in 1935.

\end{abstract}

\vyk{
\begin{keywords}
uniform sets; Borel sets; 
transfinite sequences 
\end{keywords}
}

\parf{Введение}
\las{int}

In the first half of the 1930s, 
a number of profound
and interesting results were obtained 
in descriptive set theory 
-- then a rather new branch of mathematics. 
The book \cite{L35e} by N.\,Luzin, who was 
a leader of this research area at that time, 
was devoted to the presentation and careful 
analysis of these results. 
A significant part of the book presented 
some results obtained by a young mathematician 
at that time, a student of Luzin, Petr Novikov.
In addition to those of Novikov's results, 
which, at the time of writing \cite{L35e}, 
have already been or will soon be
published in such works as 
\cite{lunoe,nov1935e,nov37ae,nov37bl},
Luzin paid attention to those studies
carried out in the doctoral thesis 
of Petr Novikov,
which were not published at that time, 
have never been published, and
are known only from their 
presentation and analysis by Luzin in \cite{L35e}.

In particular, \S\,32 and sections I--IV of 
\S\,33 of \cite{L35e} analyze well
-ordered families of uniform planar  
sets that Novikov proposed for 
<<geometric>> representation of ordinals 
$\al<\om_2$ 
(transfinite numbers of the third class and
lower, in the terminology of that time).

Considering the question of order types
of well-ordered collections of uniform    
sets, Luzin proved (Proposition~\ref{tl} below)
that these types are necessarily strictly smaller 
than $\om_2$, and poses a question (ibid., sec.\ III):
is it true that inversely, every ordinal 
$\mu<\om_2$ 
is equal to the length (= the order type) 
of some well-ordered sequence
of uniform planar sets. 

Theorem~\ref{tm} of this note
(see \S\,\ref{prel} below) answers this 
in the positive; and this is our main result. 
The uniform sets defined to prove 
the theorem, will be Borel, 
and the construction of each of them will be 
entirely effective. 
In particular,
all these uniform sets will have
G\"odel-constructible Borel codes. 
For the organization of the presentation, 
see\ at the end of \S\,\ref{prel}.

\parf{Preliminaries and the main theorem} 
\las{prel}

We consider sets on the real plane $\dR\ti\dR$, 
called simply  
{\em planar sets\/}. 
The {\em projection\/} of a planar set $P\sq \dR\ti\dR$ 
is the {\em linear\/} set 
$$
\pro P=\enx{x}{\sus y\,P(x,y)}\sq\dR, 
$$
and if $x\in\dR$ then we consider the  
{\em (vertical) section\/} 
$$
\pro P=\enx{x}{\sus y\,(\ang{x,y}\in P)}\sq\dR, 
$$
so that 
$\pro P=\enx{x}{\seq Px\ne\pu}$. 
A planar set $P$ is 
{\em uniform\/}, 
if all its sections $\seq Px$ contain 
at most one element.
This unique element will be denoted $P(x)$. 
Thus, $\seq Px=\ans{P(x)}$,  
and $P$ the graph of the function $\pro P\to\dR$. 

If $P,Q\sq\dR\ti\dR$ are uniform sets then define  
\imar{P trl Q}
$P\trl Q$ ($P$ is below $Q$) 
if the following two conditions 
are satisfied:
\ben
\nenu
\setcounter{enumi}{\value{saved}}
\itlb{uf1}%
$\pro P\sq\pro Q$ \ or \ $\pro Q\sq\pro P$;

\itlb{uf2}%
$P(x)<Q(x)$ for all $x\in \pro P\cap\pro Q$. 
\setcounter{saved}{\value{enumi}}
\een
Note that $P\trl Q$ does not determine  
which of the two inclusions is \ref{uf1}  
holds. 

Simple examples show that $\trl$ 
it is not necessarily a transitive relation.
Nevertheless we will consider collections 
of nonempty uniform planar sets
that are well-ordered  
by the relation $\trl$ --- we will call them 
{\em chains\/}, as well as those ordinals
that are their {\em lengths\/}, 
\ie\ order types.
Writing $\sis{P_\al}{\al<\mu}$, we'll mean
that this is a well-ordered chain
of uniform sets of length 
$\mu\in\Ord$, \ie\ it holds 
$P_\al\trl P_\ba$
for all $\al<\ba<\mu$.


\bpro
[\rm Novikov, Luzin \cite{L35e}, \S\,33-II]
\lam{tl}
The length\/ $\mu$ of any well-ordered chain\/ 
$\sis{P_\al}{\al<\mu}$  
of uniform set\/ $P_\al$ 
is strictly less than\/ $\om_2$. 
\epro

This result led Luzin Лузин \cite{L35e}, \S\,33-III 
to the following {\ubf reverse problem\/}: 

\begin{Qroblem}
\lam{pro}
For any ordinal\/ $\mu<\om_2$, find out, 
whether there exists a well-ordered chain of  
uniform sets of length exactly\/ $\mu$.
\end{Qroblem}

Item \ref{tm1} of 
the following theorem provides 
a positive solution to this
problem, and moreover, 
this solution
is found in the domain of {\em Borel\/}
(uniform) sets $P$,
and precisely those that satisfy $P\sq\dR\ti\dQ$, 
\ie\ only with rational ordinates.
Items \ref{tm2} and \ref{tm3} provide additional
material concerning the encoding of these Borel 
sets and efficiency of this encoding 
(in the form of the G\"odel constructibility).
In particular, claim \ref{tm3} concerns
the case when $\mu<\Oml$, where 
$\Oml$ is the first $\rL$-cardinal above 
the actual $\omi$, \ie\
if $\omi=\om_\ga^\rL$ then $\Oml=\om_{\ga+1}^\rL$. 
We also recall that $\rL$ is the class of all 
{\em constructible sets\/}. 
For concepts related to the G\"odel 
constructibility in
set theory, see, for example, 
\cite{jechmill}, \cite[\S\,2]{kl1},
or \cite[ch.\,8]{kl16e}.

\bte
[\rm the main theorem]
\lam{tm}
Suppose that\/ $\mu<\om_2$. 
Then$:$
\rm 
\sloppy
\ben
\zenu
\itlb{tm1}%
\em
there exists a chain well-ordered by\/ $\trl$, 
of length\/ $\mu$, of uniform Borel sets\/ 
$P_\al\sq\dR\ti\dQ\;\;(\al<\mu),$ 
such that in addition$:$

\rm
\itlb{tm2}%
\em
all these sets\/ $P_\al$ have constructible 
Borel codes$;$\snos
{See\ \S\,\ref{kbm} о the coding of Borel sets  
used here.}
\rm
\itlb{tm3}%
\ben
\rm
\itlb{tm3a}%
\em
if\/ $\mu<\Oml$ then the Borel codes for the sets\/ 
$P_\al$ as above can 
be chosen to form a constructible sequence, 
\rm
\itlb{tm3b}%
\em
but if\/ $\mu\ge\Oml$ then there is no any 
constructible sequence of pairwise different 
Borel codes of length\/ $\mu$.\vyk
{The restriction $\mu<\Oml$ here is necessary, 
see Lemma~\ref{avo} below. 
} 
\een
\een
\ete

It follows from \ref{tm3} that the ordinal $\Oml$ 
divides the domain of effectivity 
(understood in the sense of the existence of a 
constructible sequence of Borel codes for the 
sets in a given $\trl$-chain) and the domain 
of non-effectivity. 
The existence of such a critical ordinal in 
the semi-open interval 
$(\omi,\le\om_2]$ was predicted by Novikov;  
see Luzin \cite[\S\, 33.X]{L35e}, where the  
ordinal is denoted by $\mu^\ast$. 
Our ordinal $\Oml$ can be considered as a 
concrete realization of the Novikov--Luzin 
ordinal $\mu^*.$ 
It's role depends on the actual relationship 
between $\Oml$ and $\om_2$, of course. 
It is well known from the practice of the forcing 
method in modern set theory that
both relations $\Oml<\om_2$ and $\Oml=\om_2$ 
are compatible with the axioms of the 
Zermelo-Fraenkel theory $\ZFC$.

About the organization of the article. 
The proof of Proposition~\ref{tl} is given in 
\S\,\ref{ogr}. 
Theorem~\ref{tm} in part \ref{tm1} is 
established 
in \S\,\ref{w1} for the case of $\mu=\omi$,
and in \S\,\ref{dot} in the general case. 
Then we introduce the Borel encoding in \S\,\ref{kbm},
and prove Theorem~\ref{tm} in parts \ref{tm2} and 
\ref{tm3} in \S\,\ref{pko}.  
We discuss a corollary related to presentation of 
ordinals $\mu<\om_2$ by chains of uniform Borel 
sets in \S\,\ref{pre}, and finish with 
some concluding remarks and a problem in 
\S\,\ref{zz}.

\parf{Restriction of the length of increasing chains} 
\las{ogr}

For the convenience of the reader, 
we present here a fairly
short, though by no means obvious 
{\ubf proof of Proposition~\ref{tl}}
given by Luzin in \cite{L35e}, \S\,33--II. 
Note that Proposition~\ref{tl} {\ubf is not} 
used in the proof of Theorem \ref{tm}.

Thus let $\sis{P_\al}{\al<\mu}$ be 
a well-ordered chain of uniform sets 
$P_\al\sq\dR\ti\dR$, \ie\ $P_\al\trl P_\ba$ 
whenever $\al<\ba<\mu$. 
We have to prove that $\mu<\om_2$.

Let  
$P=\bigcup_{\al<\mu}P_\al$. 
Each section  
$\seq Px=\enx{y}{\ang{x,y}\in P}\sq\dR$ 
is well-ordered, and its order type 
$\mu_x<\omi$, 
since there are no strictly increasing  
$\omi$-sequences of reals.
Thus $\seq Px=\enx{y_{x\xi}}{\xi<\mu_x}$, 
where the numbering of the reals $y\in\seq Px$ 
is given in ascending order according to the 
usual ordering of the reals.
Let 
$
Q_\xi=\ens{\ang{x,y_{x\xi}}}{\xi<\mu_x} 
$ 
for all $\xi<\omi$, so that %
\ben
\setcounter{enumi}{\value{saved}}
\nenu
\itlb{ogr1}%
the set $Q_\xi$ are uniform, 
$Q_\xi\trl Q_\et$, \   
$\pro Q_\et\sq\pro Q_\xi$ for $\xi<\et<\omi$, \ 
and in addition 
$P=\bigcup_{\al<\mu}P_\al=\bigcup_{\xi<\omi}Q_\xi$.
\setcounter{saved}{\value{enumi}}
\een
 
We claim that
\ben
\nenu
\setcounter{enumi}{\value{saved}}
\itlb{ogr2}%
If $\xi<\omi$ then the set 
$A_\xi=\enx{\al<\mu}{Q_\xi\cap P_\al\ne\pu}$ 
is countable.
\setcounter{saved}{\value{enumi}}
\een

Suppose otherwise. 
Let $\xio<\omi$ be the least ordinal such that 
$A_\xio$ 
{\em is uncountable\/}, 
so all $A_\xi$, $\xi<\xio$, 
are countable. 
Thus by removing a countable number 
of indices from $A_\xio$, 
we get an uncountable set 
$A\sq A_\xio$ such that
\ben
\nenu
\setcounter{enumi}{\value{saved}}
\itlb{ogr3}%
if $\xi<\xio$ then $Q_\xi\cap P_\al=\pu$ for all  
$\al\in A$.
\setcounter{saved}{\value{enumi}}
\een
We assert that  
\ben
\nenu
\setcounter{enumi}{\value{saved}}
\itlb{ogr4}%
if $\al<\ba$ belong to $A$ then  
$\pro P_\al\sq\pro P_\ba$.
\setcounter{saved}{\value{enumi}}
\een

Indeed, as $\ba\in A$, there is a pair 
$\ang{x,y}\in Q_\xio\cap P_\ba$, 
and then $x\in\pro P_\ba$. 
Show that $x\nin\pro P_\al$. 
Otherwise we have $\ang{x,y'}\in P_\al$ 
for some $y'$. 
In this case, $y'<y$ is impossible, 
because then it would be 
$\ang{x,y'}\in Q_\xi$ for some $\xi<\xio$
according to \ref{ogr1}, which contradicts \ref{ogr3}. 
The equality $y'=y$ is also impossible, as 
$P_\ba\cap P_\al=\pu$ for $\al\ne\ba$. 
Finally, $y<y'$ is impossible, since it contradicts
the fact that $P_\al\trl P_\ba$ for $\al<\ba$. 
So, in fact, $x\nin\pro P_\al$. 
This implies $\pro P_\ba\not\sq\pro P_\al$,
which means that $\pro P_\al\sq\pro P_\ba$
according to \ref{uf1}, so that 
\ref{ogr4} is verified.

Continuing the proof of \ref{ogr2}, 
we take the least ordinal $\al_0\in A$ and any 
$x\in\pro{A_{\al_0}}$. 
According to \ref{ogr4}, $x\in\pro{A_{\al}}$
holds for all $\al\in A$. 
This means that there is a family of reals $y_\al$, 
$\al\in A$, for which $\ang{x,y_\al}\in P_\al$,
and at the same time $y_\al<y_\ba$ for $\al<\ba$, 
since $\al<\ba$ implies $P_\al\trl P_\ba$. 
So, we have an uncountable strictly increasing
sequence of reals $y_\al$, $\al\in A$,
in $\dR$, which is impossible. 
This contradiction completes the 
proof of \ref{ogr2}.

However, it follows from \ref{ogr2} that this chain 
$\sis{P_\al}{\al<\mu}$ contains no more than 
$\ali$ sets, which means $\mu<\om_2$. 
This ends the proof of Proposition \ref{tl}.

\parf{The main theorem for the first uncountable 
ordinal} 
\las{w1}

Here we present the
{\ubf proof of Theorem~\ref{tm}},
only in part  
\ref{tm1}, {\ubf for the case of $\mu=\omi$}. 
This result will be an important step for
the general case. 

As we consider mainly
Borel sets $B\sq\dR\ti\dQ$, 
the {\em decomposition\/} of such 
sets will be used, 
onto {\em horizontal sections\/}, \ie, 
also obviously Borel sets:
\setcounter{equation}{\value{saved}}
\bul
{Wq}
{\textstyle
\hors Br=\enx{x}{\ang{x,r}\in B}\sq\dR, \;
\text{ where }\,r\in\dQ, 
\text{  so }
B=\bigcup_{r\in\dQ}\big(\hors Br\ti\ans r\big).
}

We begin with a standard lemma. 
As usual, $\dQ=$ rational numbers.

\ble
\lam{dotL1}
There is such a Borel set\/ 
$G\sq\dR\ti\dQ$ 
that for every \/$Z\sq\dQ$ there exists 
a real\/ $x\in\dR$ satisfying\/ 
$Z=\seq Gx:=\enx{q\in \dQ}{\ang{x,q}\in G}$.
\ele
\bqf
We fix a recursive enumeration of rational
numbers, $\dQ=\enx{r_n}{n<\om}$. 
For $x\in\dR$, let $A_x\sq\om$ be the set
of all $n$ such that the decomposition of 
the fractional part of $x$ into
a binary fraction 
has $0$ at position $2n+1$. 
The set $G=\enx{\ang{x,r_n}}{n\in A_x}$ 
is as required.
\eqf

We fix such a set of $G$ for further consideration. 

You may notice that $G$ is one of the options. 
{\em of the Lebesgue binary sieve\/} 
\cite{leb1905},
for which see \eg\ \cite{L35e}, \S\,17.
Define, for each ordinal $\xi<\omi$,
\bul
{UD}
{\left.
\bay{lcl}
D_\xi &=& \{x\in \dR:
\text{the section 
$\seq Gx$ 
has a well-ordered} \\[0.0ex] 
&&\phantom{\{x\in \dR:{}}
\text{initial segment of length ${>}\,\xi$}\};
\\[0.9ex]
U_{\xi} &=& \enx{\ang{x,r}\in D_{\xi}\ti\dQ}
{r
\text{ is the $\xi$th largest real in }
\seq Cx};
\eay
\right\}
}
and finally $\cD = \ens{D_\xi}{\xi<\omi}$.

The next lemma implies Theorem~\ref{tm}, 
claim \ref{tm1}, in case $\mu=\omi$:

\ble
\lam{dotB}
Sets\/ $U_\xi\sq\dR\ti\dQ$, 
$\xi<\omi$, are uniform and
form a \/$\trl$-increasing chain\/ 
$\sis{U_\xi}{\xi<\omi}$   
of lengths\/$\omi\,;$ also\/
$D_\xi=\pro U_\xi$
and\/ 
$D_\et\sq D_\xi$ for \/$\xi<\et$. 

All sets \/$U_\xi$ and \/$D_\xi$ are  
nonempty and Borel.
\ele

\bqf
Only the Borelness needs to be proved, 
the rest is obvious.
We prove that $U_\xi$ is Borel by induction 
on $\xi$. 
Assume for simplicity that the variables $p,q,r$
denote only rational numbers.
For $\xi=0$:
\vyk{
\bul
{bo1}
{
\ang{x,q}\in U_0\,\leqv\,
\ang{x,q}\in C \land 
\kaz q'<q\,
(
\ang{x,q'}\nin C
), 
}
}
\bul
{bo1}
{
x\in \hors{U_0}r\,\leqv\,
x\in \hors{G}r  \land 
\kaz q<r\,
\big(
x\nin \hors{G}q
\big), 
}
which implies the Borelness of all sections  
$\hors{U_0}r$
and $U_0$ itself as well, since
the set $G$ is Borel
by lemma \ref{dotL1}, and the quantifier $\kaz q$
on the right-hand side 
is limited to the countable set $\dQ$.

Now assume that $\xi>0$, and all sets $U_\et$, 
$\et<\xi$ are Borel. 
Then
\bul
{bo2}
{
\left.
\bay{l}
x\in \hors{U_\xi}r\;\leqv\;
x\in \hors{G}r \,\land\, 
\kaz\et<\xi \,\sus q<r\, 
\big(
x\in \hors{U_\et}q
\,\land {} \\[0.5ex]
\hspace*{30ex}
{}\land \,\kaz p\,(q'<p<q\imp x\nin \hors Gp)
\big) 
\eay
\right\}
}
in case $\xi=\et+1$, 
while if $\xi$ is a limit ordinal then
\bul
{bo3}
{
\left.
\bay{l}
x\in \hors{U_\xi}r\;\leqv\;
x\in \hors Gr\, \land\, 
\kaz\et<\xi\,\sus q<r\, 
\big(
x\in \hors{U_\et}r
\big) 
\land{}\\[0.5ex]
\hspace*{20ex}
{}\land\, 
\kaz q<r\,\sus\et<\xi\,
\big(
x\in \hors Gq\imp x\in \hors{U_\et}q 
\big). 
\eay
\right\}
}%
In both cases $U_\xi$ is Borel since such are 
both $G$ and all sets 
$U_\et$, $\et<\xi$. 

Finally, each set $D_\xi=\pro U_\xi$ 
is Borel because
\bul
{bo44}
{%
x\in D_\xi\,\leqv\,
\sus r\in\dQ\,\big(x\in \hors{U_\xi}r\big),
\setcounter{saved}{\value{equation}}
}%
which completes the proof of Lemma \ref{dotB}.
\eqf

\parf{Proof of the main theorem, part \ref{tm1}, 
in the general
case} 
\las{dot}

We begin with an important special case 
of {\em transitivity\/} of the relation $\trl$.

\ble
\lam{tran}
Assume that\/ 
$P\,{\trl}\hspace*{0.8ex}U\hspace*{0.2ex}
{\trl}\hspace*{0.7ex}Q$ are uniform sets. 
Then the conjunction of\/ \ref{tran2} below
and\/ \ref{uf1} above in\/ {\rm\S\,\ref{prel}} 
suffices 
for\/ $P\trl Q$ to hold$:$\vim
\ben
\nenu
\setcounter{enumi}{\value{saved}}
\itlb{tran2}
$\pro P\sq\pro U$ 
and/or\/ $\pro Q\sq\pro U$. 
\setcounter{saved}{\value{enumi}}
\een
\ele

\begin{qroof}
To prove \ref{uf2} take any 
$x\in\pro P\cap\pro Q$ by \ref{uf1}. 
Then $x\in \pro U$ by \ref{tran2}, 
and we are done.  
\end{qroof}

Say that a chain of uniform sets $P_\al$
has the {\em property of $\cD$-projections} if 
$\pro P_\al\in\cD$ for all $\al$. 
The chain 
$\sis{U_\xi}{\xi<\omi}$ obviously has this property.

\ble
\lam{dotL3}
Let\/ $\mu<\om_2$. 
There is a\/ $\trl$-increasing chain of
uniform Borel sets\/ $P\sq\dR\ti\dQ$,
of length $\mu$, 
with the property of\/ $\cD$-projections.
\ele

\bqf
According to lemma \ref{dotB}, 
the chain  
$\sis{U_\xi}{\xi<\omi}$ handles the case 
$\mu\le\omi$. 
Next, we argue by induction.\vom 

{\ubf Beginning of the inductive step.}
Fix an ordinal $\mu;\;\omi<\mu<\om_2$. 
\ben
\setcounter{enumi}{\value{saved}}
\nenu
\itlb{fp}%
Fix an enumeration 
$\mu=\enx{\nu_\xi}{\xi<\omi}$ of all smaller ordinals 
$\nu<\mu$. 
(We will return to the analysis of this action below.)
\setcounter{saved}{\value{enumi}}
\een
For each ordinal $\nu_\xi$,
by the inductive hypothesis,
there is
\ben
\setcounter{enumi}{\value{saved}}
\nenu
\itlb{FF}%
a $\trl$-increasing chain 
$\sis{F^\xi_\al}{\al<\nu_\xi}$, 
of length $\nu_\xi$, 
of uniform Borel sets  
$F^\xi_\al\sq\dR\ti\dQ$, 
with the property of $\cD$-projections. 
\setcounter{saved}{\value{enumi}}
\een
Next, we're going to insert each chain 
$\sis{F^\xi_\al}{\al<\nu_\xi}$  
between the uniform sets $U_\xi$ and $U_{\xi+1}$
by a vertical shift. 
To this end, we first define 
\pagebreak[0]
\setcounter{equation}{\value{saved}}
\bul
{defQ}
{%
\left.
\bay{l}
Q^\xi_\al=F^\xi_\al\res D_{\xi+1}:=
P^\xi_\al\cap(D_{\xi+1}\ti\dQ)
\sq\dR\ti\dQ\,,\\[0.4ex]
\text{or equivalently,}\quad 
\hors{Q_\al^\xi}r=\hors{F_\al^\xi}r\cap D_\xip
\,,\;\kaz r\,,
\eay
\right\}
} 
for all $\al<\nu_\xi$. 
As $D_{\xi+1}\in\cD$, and the chain 
$\sis{F^\xi_\al}{\al<\nu_\xi}$
has the property of $\cD$-projections, all reduced
sets $Q^\xi_\al$ are nonempty and form 
$\trl$-ascending chain 
$\sis{Q^\xi_\al}{\al<\nu_\xi}$
of the same length $\nu_\xi$ and
with the property of $\cD$ projections, and 
$\pro Q^\xi_\al\sq D_{\xi+1}\sq D_\xi,\:\kaz \al.$

All sets $Q^\xi_\al$ are Borel, since such are 
all $F^\xi_\al$ 
and $D_\xi$ (by Lemma \ref{dotB}). 

Let $\xi<\omi$ and $x\in D_{\xi+1}$. 
By Lemma \ref{dotB}, there are  
unique pairs 
$\ang{x,u^x_\xi}\in U_\xi$ and 
$\ang{x,u^x_{\xip}}\in U_{\xip}$,
moreover, $u^x_{\xi}<u^x_{\xip}$ are rational.

For any pair  $u<v$ of rational numbers,
define an {\em order-preserving\/}  
bijection 
$H[u,v]:(u,v)\ontor\dR$ 
of an open interval $(u,v)$
onto $\dR$:
\bul
{huv}
{H[u,v](y)=
\left\{
\bay{rcl}
\fras{y-c}{v-y}
&\text{in case}&
c\le y<v \;; \\[1.7ex]
\fras{y-c}{u-y}
&\text{in case}&
u< y\le c\;; 
\eay
\right|
\text{ where \ }\textstyle c=\fras{u+v}2\,.
}
Clearly $H[u,v]$ preserves the rationality, 
\ie\ $r\in\dQ\eqv H[u,v](r)\in\dQ$.

Let again $\xi<\omi$. 
If $\al<\nu_\xi$, then the uniform
set $Q^\xi_\al\sq\dR\ti\dQ$  
satisfies $\pro Q^\xi_\al\sq D_{\xi+1}$ 
by the above. 
Define the {\em vertical shift\/}
\bul
{defS}
{S^\xi_\al=\enx{\ang{x,r}\in D_{\xip}\ti\dQ}
{ u^x_{\xi}<r<u^x_{\xip}\land 
\ang{x,H[u^x_{\xi},u^x_{\xip}](r)}\in Q^\xi_\al}
}%
of $Q^\xi_\al$ into the area between $U_\xi$ and
$U_{\xip}$. 
As the abscissas do not change, 
and the ordinates
maintain order and rationality, we conclude   
that the sets $S^\xi_\al\sq D_{\xip}\ti\dQ$, 
$\al<\nu_\xi$, are uniform
and form a $\trl$-increasing chain 
$\sis{S^\xi_\al}{\al<\nu_\xi}$ of length $\nu_\xi$ and  
with the property of $\cD$-projections, and 
in addition
\bul
{mez}
{\pro S^\xi_\al=\pro Q^\xi_\al\sq D_{\xip}\sq D_\xi  
\qand U_\xi\trl S^\xi_\al\trl U_{\xip}
\setcounter{saved}{\value{equation}}
}
for all $\al<\nu_\xi$ by construction. 

Prove that these new sets $S^\xi_\al$ are Borel. 
%
%
Indeed by construction, 
\bul
{sxa}
{\hors{S^\xi_\al}{r}=
\bigcup_{u,v,w\in\dQ\land u<v\land w=H[u,v](r)}
\hors{U_\xi}{u}\cap\hors{U_\xip}{v}
\cap\hors{Q^\xi_\al}{w}.
\setcounter{saved}{\value{equation}}
}
for any $r\in\dQ$. 
Yet the sets
${U_\xi},\,{U_\xip},\,{Q^\xi_\al}\sq\dR\ti\dQ$ are  
Borel, therefore their horizontal sections
$\hors{U_\xi}{u},\,\hors{U_\xip}{v},\,
\hors{Q^\xi_\al}{w}$
are Borel as well. 
Therefore, according to \eqref{sxa}, the sections 
$\hors{S^\xi_\al}{r}$ are Borel, too. 
This implies that the sets $S^\xi_\al$ themselves 
are Borel by \eqref{Wq}, as required.

Let's now analyze the whole family
\pagebreak[0]
\bul
{ww}
{\smu_\mu=\enx{U_\xi}{\xi<\omi}\cup
\enx{S^\xi_\al}{\xi<\omi\land\al<\nu_\xi} 
\setcounter{saved}{\value{equation}}
}
of uniform subsets in $\dR\ti\dQ$,
considered for the inductive step $\mu$.
Prove that  
\ben
\setcounter{enumi}{\value{saved}}
\nenu
\itlb{21}%
$\smu_\mu$ is a $\trl$-chain of length 
$\mu'=\sum_{\xi<\omi}(1+\nu_\xi)\ge\mu$. 
\setcounter{saved}{\value{enumi}}
\een

{\em Fact 1\/}. 
We already know that both $\sis{U_\xi}{\xi<\omi}$ 
and 
$\sis{S^\xi_\al}{\al<\nu_\xi}$ 
for any $\xi<\omi$, are 
$\trl$-increasinging chains of uniform 
sets of corresponding lengths, 
and with the property of $\cD$-projections, 
and \eqref{mez} holds.\vom

{\em Fact 2\/}. 
If $P,Q\in\smu_\mu$ then the projections 
$\pro P$ and $\pro Q$ 
belong to $\cD$, hence  
$\pro P\sq\pro Q$ or $\pro Q\sq\pro P$ by
\ref{uf1} of the definition of $\trl$. \vom

{\em Fact 3\/}. 
If $\xi<\et<\omi$ and $\al<\nu_\xi$ then 
$S^\xi_\al\trl U_{\xip}\trl U_\et$ by Fact~1,
$\pro U_\et\sq\pro{U_{\xi+1}}$ by Lemma \ref{dotB},
thus $S^\xi_\al\trl U_\et$ by Lemma \ref{tran}
and Fact~2.\vom

{\em Fact 4\/}. 
If $\xi<\et<\omi$ and $\ba<\nu_\et$ then 
$U_\xi\trl U_\et\trl S^\et_\ba$ by Fact 1,
and   
$\pro S^\et_\ba\sq\pro{U_{\et}}$
by \eqref{mez},
so that $S^\xi_\al\trl U_\et$ by Lemma~\ref{tran}
and Fact 2.\vom

{\em Fact 5\/}. 
If $\xi<\et<\omi$ and $\al<\nu_\xi$, $\ba<\nu_\et$,
then $S^\xi_\al\trl S^\et_\ba$. 
Here we have  
$S^\xi_\al\trl U_\et\trl S^\et_\ba$
by Fact 3,
and  
$\pro S^\et_\ba\sq\pro{U_{\et}}$ by \eqref{mez},
whence $S^\xi_\al\trl S^\et_\ba$ holds 
again by Lemma~\ref{tran} and Fact 2.\vom

As a consequence of Facts 1--5, we get \ref{21} 
for the family \eqref{ww},
\ie
\setcounter{equation}{\value{saved}}
\bul
{wwp}
{\smu_\mu=\enx{P_\al}{\al<\mu'} 
\setcounter{saved}{\value{equation}}
}
in ascending order of $\trl$.
However, $\mu'\ge\mu$ is obvious, 
so that a chain of length
exactly $\mu$ is obtained by simply cutting the 
<<tail>> $\smu_\mu$ in \eqref{wwp} from   
$\mu$ to $\mu'$. 
This {\ubf completes the inductive step} and the proof
of Lemma~\ref{dotL3}.
\eqf

Finally, Lemma~\ref{dotL3} obviously {\ubf implies
Theorem ~\ref{tm}, in part \ref{tm1}}. \ \ \ \qed

\parf{Encoding Borel sets} 
\las{kbm}

Here we recall the basic concepts in connection 
with Borel codes, which appear in
the formulations and will be used
in {\ubf the proofs of the statements \ref{tm2}
and \ref{tm3} of Theorem~\ref{tm}} 
in the next section. 
We consider  the set $\omi\lom$ of all 
tuples (finite sequences) of countable
ordinals. 
Further:
\bit
\item[$-$] 
$s\su t$ means that the tuple $t$ 
is a proper extension of the tuple $s$, 
 
\item[$-$] 
 $\puk$ is the empty tuple, 
$\ang{\al_1,\dots,\al_n}$ is the tuple with 
terms 
$\al_1,\dots,\al_n$,  

\item[$-$] 
$s\we \al$ is obtained by adjoining  
the rightmost term $\al$ to the tuple $s$,

\item[$-$] 
a set $T\sq\omi\lom$ is a 
{\em tree\/}, 
if $s\su t\in T\imp s\in T$, 

\item[$-$] 
$\Max T$ is the set of all 
{\em endpoints\/} of a given tree $T$, 

\item[$-$] 
a tree $T$ is {\em \wef\/}, if it has no  
{\em infinite branches\/}, \ie\ there is no 
function $b:\om\to\omi$ 
such that $b\res n\in T$ for all $n$. 
\eit

Finally, let 
{\em a Borel code\/} (for the space $\dR$),
be any pair $\ang{T,d}$, where 
$\pu\ne T\sq\omi\lom$ -- 
{\em is a finite or countable\/}
\wf\ tree, and $d\sq T\ti\dQ\ti\dQ$.
In this case, the \wf ness of the tree $T$ 
allows to uniquely define
the Borel set $[T,d,s]\sq\dR$ 
for each $s\in T$ so that:
\ben
\Renu
\itlb{b1}\msur%
$[T,d,s]=
\dR\bez \bigcup_{\ang{s,p,q}\in d} (p,q),$  
in case $s\in\Max T$;%
 
\itlb{b2}\msur%
$[T,d,s]= 
\dR\bez\bigcup_{s\we\al\in T}[T,d,s\we\al]$, 
in case $s\in T\bez\Max T$; 

\itlb{b3}%
finally $[T,d]=[T,d,\puk]$; 
\een
where, as usual, $(p,q)=\enx{x\in\dR}{p<x<q}$
in \ref{b1} is a rational open
interval of the real line $\dR$,
empty for $p\ge q$.

Thus the scheme \ref{b1}, \ref{b2}, \ref{b3}
defines the set $[T,d]\sq\dR$ 
from rational intervals,
by the operation of the complement 
to the countable union, \ie\ a Borel set. 
Conversely every Borel set $X\sq\dR$ 
admits a Borel code $\ang{T,d}$ 
(with a countable tree $T$!) for which $X=[T,d]$.%
\snos
{\label{komm}%
This definition 
corresponds to the topology of the real line 
$\dR.$
For encoding Borel sets, for example, 
{\em of the Baire space\/} $\dI=\om^\om$ 
we should take the sets $d\sq T\ti\om\lom,$
and also change \ref{b1} to the form 
$[T,d,s]=
\dI\bez \bigcup_{\ang{s,\sg}\in d}I_\sg,$ 
where $I_\sg=\enx{a\in\dI}{\sg\su a}$.} 

As for the encoding of {\em planar\/}
Borel sets,
fortunately, we do not need to consider 
this question in all generality, 
since in fact only planar sets  
$U\sq\dR\ti\dQ$ occur  
in the proof of Lemma~\ref{dotL3} and
Theorem~\ref{tm}. 
Let a  
{\em Borel multicode\/} be any indexed system 
of Borel codes 
$c=\sis{{T_r,d_r}}{r\in\dQ}$.
Then we define
\bce
$[c]=
\bigcup_{r\in\dQ}[T_r,d_r]\ti\ans r
\sq\dR\ti\dQ$,

\ece%
so that $[c]=U$ in case $\hors Ur=[T_r,d_r]$ 
for all $r\in\dQ$.

\vyk{
\bdf
\lam{bmc}
$\BC$ is the set of all 
Borel codes; $\BMC$ is the set of all 
Borel multicodes.
\edf
}

With this definition,  
Borel codes and multicodes are 
\rit{hereditarily countable\/} sets,
which plays an essential role in some
definability issues. 
If we allowed rational intervals per se  
instead of pairs of their endpoints in the 
conditions for $d$
(which would seem to be a simpler and more natural
solution),
then this hereditary countability would 
disappear, of course.


\parf{Defining codes to prove the main
theorem} 
\las{pko}

Having outlined these standard definitions 
related to Borel codes, 
let us return to Theorem ~\ref{tm}. 
The purpose of this section is  
{\ubf the proof of statements 
\ref{tm2} and \ref{tm3} of the theorem} 
for Borel sets \eqref{ww}/\eqref{wwp} 
constructed above in \S\,\ref{dot} for the 
proof of the theorem in part \ref{tm1}.
The next intermediate
result gives Borel codes for 
sets $U_\xi$ and $D_\xi$ from \eqref{UD},
which belong 
to the class $\rL$ of 
{\em G\"odel constructible\/} sets. 

\ble
\lam{cor1}
It is true for the sets\/ 
$U_\xi,\,D_\xi \,\;(\xi<\omi)$ in\/ 
\eqref{UD} that\/$:$\vim 
\ben
\renu
\itlb{U1}%
there exist Borel multicodes\/  
$c_\xi=\sis{{T^\xi_r,d_r^\xi}}{r\in\dQ}\in\rL$
for\/ $U_\xi$, and
Borel codes\/  
$\ang{T'_\xi,d'_\xi}\in\rL$ for\/ $D_\xi$, 
such that$:$

\itlb{U2}%
the\/ $\omi$-sequences of these codes are 
constructible as well.
\een
\ele
\bqf[\rm sketch]
Equations \eqref{bo1}--\eqref{bo3}  
allow to effectively
define, by transfinite induction,  
Borel multicodes 
$c_\xi=\sis{{T^\xi_r,d_r^\xi}}{r\in\dQ}$
for sets $U_\xi$, 
\ie\ $U_\xi=[c_\xi]$,
and then, using \eqref{bo44},   
Borel codes $\ang{T'_\xi,d'_\xi}$  
for sets $D_\xi=\pro U_\xi$ as well. 
We should start by defining a constructible
multicode for the initial set  
$G$ given by the proof of Lemma \ref{dotL1}; 
we'll leave this as a simple exercise.
Inductive construction of all 
these codes
is absolute  for the class $\rL$ of all 
constructible sets, 
\ie\ gives the same result in $\rL$
and in the universe of all sets.
\eqf

Lemma \ref{cor1} together with lemma~\ref{dotB} above
give a complete proof of Theorem \ref{tm} with
all three of its parts \ref{tm1}, \ref{tm2}, 
\ref{tm3} for $\mu=\omi$. 
{\ubf 
The proof for the general case $\mu<\om_2$}
is based on the next lemma,
similar to Lemma \ref{cor1}\ref{U1}, but relevant 
to sequences of Borel sets
from the proof of Lemma~\ref{dotL3}.

\ble
\lam{cor2}
Let\/ $\omi<\mu<\om_2$. 
In the context of the notation and assumptions
of the inductive step   
in the proof of Lemma~\ref{dotL3}, suppose 
that
\ben
\fenu
\itlb{2*}%
there exist
Borel multicodes\/ 
$\vpi^\xi_\al \in\rL$ 
for the sets\/ $F^\xi_\al$ as in \ref{FF}.
\een
Then there exist Borel multicodes 
$\pi_\al\in\rL$
for the resulting sets\/ $P_\al$ in  
\eqref{wwp} in   
the proof of Lemma~\ref{dotL3}.\vyk
{\lam{nev}%
A statement like \ref{U2} of Lemma \ref{cor1} 
is not possible here, see Lemma~\ref{avo} below. 
}
\ele

\bqf[\rm sketch]
First, we define intermediate multicodes  
$\qoppa^\xi_\al\in\rL$
for sets $Q^\xi_\al$ from codes \ref{2*} 
using relations \eqref{defQ} in \S\,\ref{dot},
then multicodes 
$\sg^\xi_\al\in\rL$
for sets \/$S^\xi_\al$ using relations 
\eqref{sxa} in the same place, which, 
along with the multicodes
for the sets $U_\xi$, provided  
by Lemma~\ref{cor1}, become the required 
multicodes $\pi_\al\in\rL$ for sets $P_\al$
after renumbering during the transition 
from \eqref{ww} to \eqref{wwp}.
Once again, the constructibility of all these 
multicodes follows from the obvious absoluteness.
\eqf

This completes {\ubf the proof of   
Theorem~\ref{tm}, claim \ref{tm2}}. \ \ \ \vtm

Turning to {\ubf part \ref{tm3} of the theorem},
note that the constructions given in 
\S\S\,\ref{w1} and \ref{dot} contain only one 
``ineffective'' action involving  
an arbitrary choice, namely, 
the choice of a specific
enumeration $\mu=\enx{\nu_\xi}{\xi<\omi}$ of all
ordinals $\nu<\mu$ in the agreement \ref{fp} in 
the proof of Lemma~\ref{dotL3}. 
This, of course, does not allow to strengthen
Lemma~\ref{cor2} by requiring the constructibility 
of the resulting $\mu$-sequence of codes for 
sets \eqref{ww}/\eqref{wwp} in the proof
of Lemma~\ref{dotL3}.

Fortunately, this does not affect 
the constructibility of the 
codes for sets \eqref{ww}/\eqref{wwp},
since, for $\xi<\omi$ fixed, the formulas 
\eqref{defS} and \eqref{sxa} depend only on the value 
$\nu_\xi$ for this $\xi$, and not on the entire
sequence $\sis{\nu_\xi}{\xi<\omi}$.

As for the case when $\mu<\Oml$,
considered in statement \ref{tm3} 
of Theorem~\ref{tm}, the inequality $\mu<\Oml$
allows makes it possible to select a \rit{specific}  
enumeration $\mu=\enx{\nu_\xi}{\xi<\omi}$,
namely, the \rit{smallest one\/}  
of all such enumerations, 
in the sense of the canonical well-ordering 
$<_{\rL}$ of the constructible universe $\rL$. 
And then the construction of a sequence of sets 
\eqref{ww}/\eqref{wwp}, and the according codes   
$\pi_\al$ in the proof of Lemma \ref{cor2} 
becomes fully effective and absolute for $\rL$.\vom

This completes {\ubf the proof
of Theorem~\ref{tm} in
its part \ref{tm3a}}.\vom\  

We finish this section with {\ubf the proof
of Theorem~\ref{tm} in its part \ref{tm3b}}. 


Consider the set $\BMC$ of all Borel multicodes 
as defined in \S\,\ref{kbm}. 
Suppose that $\omi=\om_\ga^\rL$; 
then by definition $\Oml=\om_{\ga+1}^\rL$. 
On the other hand, the set $\BMC\cap\rL$ 
of all {\em constructible\/} multicodes 
(this is different from $(\BMC)^\rL$!) 
has cardinality $\aleph_\ga^\rL$ in $\rL$ 
since GCH holds in $\rL$. 
Therefore no {\em constructible\/} $\Oml$-sequence 
of pairwise different Borel multicodes 
in $\BMC\cap\rL$ can exist.\vom

This completes {\ubf the proof
of Theorem~\ref{tm} in
its last part \ref{tm3}}.\vtm\ \ \ \qed

\parf{On the effective representation of ordinals} 
\las{pre}
In this short section, we will point out one rather
fundamental, albeit somewhat vague
consequence of Theorem\ref{tm}, 
\ie\ rather a consequence of the
effectiveness of the 
construction of a $\trl$-increasing sequence
of uniform Borel sets
of a given length $\mu<\om_2$,
in the course of the proof
of Theorem~\ref{tm}.

It is clear that the von Neumann ordinals denote   
transfinite increase, so to speak, vertically 
in the hierarchy of the set-theoretic
universe, while the real numbers, their
sets, for example, Borel, and then for example 
transfinite sequences of these sets
symbolize only the expansion of the universe at
several initial steps of the ``vertical'' hierarchy. 
Therefore, it is quite natural for the foundations of mathematics, that the question arises about 
the representation, modeling, or, if you prefer, 
``naming''
ordinals of a particular magnitude by means 
of objects related to the real line. 

For example, the set $G\sq\dR\ti\dQ$, \ie\ 
the binary Lebesgue sieve, given
by Lemma~\ref{dotL1}, 
defines a representation of countable
ordinals (domain $\xi<\omi$), in which   
any given countable ordinal $\xi$ is represented 
by the Borel set $U_\xi\sq\dR\ti\dQ$, 
defined by the formula \eqref{UD}.
By the way, these sets are nonempty and 
pairwise disjoint.
Another representation can be given by Borel 
{\em linear\/} sets 
$D_\xi=\pro U_\xi\sq\dR$. 
However, the sets $D_\xi$, although nonempty, 
are not disjoint, but on the contrary, 
they are nested in one another. 
This can be fixed by taking the differences 
$E_\xi=D_\xip\bez D_\xi$, which are also
Borel sets, non-empty and pairwise disjoint.

Thus, there is an effective
representation (in different but related versions)
of countable ordinals by Borel sets, whose origins 
can be traced to the old work of 
Lebesgue~\cite{leb1905}.

Our Theorem~\ref{dotL3} in part \ref{tm3} 
yields an effective
representation in a much wider area $\mu<\Oml$. 
Namely, an effective  
representation of the ordinal $\mu<\Oml$ is 
given by the  
$\trl$-increasing $\mu$-sequence
of Borel uniform sets, and a constructive
sequence of codes for them,
the existence of which is given for this $\mu$ 
by Theorem~\ref{tm}\ref{tm3}.

As for the ordinals $\mu$ in the interval 
$\Oml\le\mu<\om_2$ 
(provided that strictly $\Oml<\om_2$),
the most reasonable efficient code 
for such a $\mu$ in this context
seems to be 
the set of all 
$\mu$-sequences of constructive
multicodes for $\trl$-chains of length $\mu$
given by Theorem~\ref{tm} in part \ref{tm1}.

\parf{Concluding remarks} 
\las{zz}

Our Theorem~\ref{tm} closes
the long-known classical problem
of Petr Novikov and Luzin on
the lengths of transfinite sequences
of uniform planar sets, and in the strongest 
form of {\em Borel\/} uniform sets of
the space $\dR\ti\dQ$  
(\ie\ with rational ordinates) 
with constructive Borel codes. 
We expect that the results obtained will find
applications in modern research in
descriptive set theory.

We also expect that our methods of effective
transfinite constructions will make a definite
contribution to the modern theory of generalized
computability on uncountable
structures and the theory of information 
transmission between
structures on ordinals and Borel
structures associated with the real line, 
as in recent works \cite{ham_effmat,hamTM}.

We finish with a problem. 
Coming back to \ref{tm3} of Theorem~\ref{tm}, the 
ordinal $\Oml$ can be seen as a separator between 
the lengths of effectively existing and effectively 
non-existing $\trl$-chains of uniform Borel sets, 
with the effectivity being understood in terms of the 
existence of a 
related constructible sequence of Borel codes. 
But there are other measures of effectivity, like 
\eg\ the ordinal-definability. 
The following problem arises:

\begin{Qroblem}
Find a model in which\/ $\Oml<\om_2$, 
and there is no {\em ordinal-definable\/} 
(or even {\em real-ordinal-definable\/}) chains of  
Borel uniform sets of length\/ $\ge{\Oml}$.
\end{Qroblem}

\begin{ackn}
The authors are grateful to Mirna D\v zamonja for 
an interesting
discussion and valuable comments.
\end{ackn}

\let\section\subsection

\renek{\refname} {References}

\bibliographystyle{plain}

\small



\end{document}